# ON RATES OF CONVERGENCE FOR POSTERIOR DISTRIBUTIONS IN INFINITE-DIMENSIONAL MODELS

By Stephen G. Walker,[1] Antonio Lijoi[2] and Igor Prünster[3]

*University of Kent, Università degli Studi di Pavia
and Università degli Studi di Torino*

This paper introduces a new approach to the study of rates of convergence for posterior distributions. It is a natural extension of a recent approach to the study of Bayesian consistency. In particular, we improve on current rates of convergence for models including the mixture of Dirichlet process model and the random Bernstein polynomial model.

**1. Introduction.** Recently, there have been many contributions to the theory of Bayesian consistency for infinite-dimensional models. Most of these adopt the "frequentist" (or "what if") approach, which consists of generating independent data from a "true" fixed density $f_0$ and checking whether the sequence of posterior distributions accumulates in Hellinger neighborhoods of $f_0$. The determination of sufficient conditions for Hellinger consistency has been the main goal of a number of recent papers such as, for example, [1, 2, 5] and [12]. A summary is provided in [8]. Their results rely upon the use of uniformly consistent tests, combined with the construction of suitable sieves and computation of metric entropies. An alternative method for solving the problem can be found in [14], where a sufficient condition in terms of the summability of prior probabilities is provided.

Here, we consider the allied problem of determining rates of convergence, that is, the determination of a sequence $(\varepsilon_n)_{n\geq 1}$ such that $\varepsilon_n \downarrow 0$ and

$$\Pi_n(\{f : d(f, f_0) > M\varepsilon_n\}) \to 0$$

Received October 2004; revised September 2006.
[1]Supported by an EPSRC Advanced Research Fellowship.
[2]Also affiliated with CNR-IMATI, Milan, Italy. Supported in part by MUR grant 2006/134525.
[3]Also affiliated with Collegio Carlo Alberto and ICER, Turin, Italy. Supported in part by MUR grant 2006/133449.

for any constant $M > 0$. The above-displayed convergence can be understood either as convergence in $F_0^\infty$-probability or as almost sure-$F_0^\infty$, where $F_0$ denotes the probability distribution associated with $f_0$ and $F_0^\infty$ is the infinite product distribution. Among recent papers dealing with this topic, we mention [4, 6, 7] and [13]. The key to these papers is the construction of a sieve and the use of entropies. The ultimate rate of convergence achieved depends on two quantities: the concentration rate, which depends on the prior mass assigned to suitable neighborhoods of $f_0$, and the growth rate of the Hellinger entropy. A recent contribution, relying upon information theory, is given in [15]. The aim of the present paper is to tackle the problem based on the approach of Walker [14], which leads to improvements in the examples we consider.

In Section 2, we first derive a useful bound for the posterior probability on the sets of interest and then prove a general theorem for the determination of rates, which relies upon two conditions. In Section 3, the normal mixture of Dirichlet process and random Bernstein polynomials are considered and currently known rates are improved.

**2. Posterior convergence rates.** Consider a sequence of observations $(X_n)_{n \geq 1}$, each taking values in some Polish space $\mathbb{X}$ endowed with the Borel $\sigma$-algebra $\mathscr{X}$. If $\mathbb{F}$ indicates the space of probability density functions with respect to some $\sigma$-finite measure $\lambda$ on $\mathbb{X}$, then the Hellinger metric $h$ on $\mathbb{F}$ is defined by

$$h(f,g) = \left\{ \int_{\mathbb{X}} (\sqrt{f(x)} - \sqrt{g(x)})^2 \, \lambda(dx) \right\}^{1/2}$$

for any $f$ and $g$ in $\mathbb{F}$, and we set $\mathscr{F}$ to be the Borel $\sigma$-algebra of $\mathbb{F}$. Suppose that $\Pi$ stands for a prior distribution on $(\mathbb{F}, \mathscr{F})$. Then the posterior distribution, given the observations $(X_1, \ldots, X_n)$, coincides with

$$\Pi_n(B) = \frac{\int_B \prod_{i=1}^n f(X_i) \Pi(df)}{\int_{\mathbb{F}} \prod_{i=1}^n f(X_i) \Pi(df)}$$

for all $B$ in $\mathscr{F}$. We assume that there exists a "true" density function $f_0$ such that the $X_n$'s are i.i.d. from $f_0$. A sequence of posterior distributions $\Pi_n$ is said to be *Hellinger consistent* at $f_0$ if the posterior mass on sets of the type $A_\varepsilon := \{f : h(f, f_0) > \varepsilon\}$ becomes negligible as the sample size $n$ increases. The approach introduced in [14] relies upon the construction of a suitable covering of $A_\varepsilon$ by Hellinger balls of radius $\phi < \varepsilon$. The prior mass on these balls must be such that the sum of their square roots is finite. This entails consistency. Then, when dealing with rates, it is natural to refine the set $A_\varepsilon$ to $A_{\varepsilon_n} = \{f : h(f, f_0) > \varepsilon_n\}$ and to consider a covering $\{A_{n,j} : j = 1, 2, \ldots\}$ of $A_{\varepsilon_n}$, where each $A_{n,j}$ has radius $\phi_n \in (0, \varepsilon_n)$. Consequently, we now define

$$K_{\varepsilon_n} = \sum_{j \geq 1} \Pi(A_{n,j})^{1/2},$$

a quantity that will be heavily relied on in this paper.

Before stating the preliminary result, let us introduce some notation. Let $L_{0,j}^{(n)} = \sqrt{\Pi(A_{n,j})}$ and for any $k \geq 1$, let

$$L_{k,j}^{(n)} := \sqrt{\int_{A_{n,j}} R_k(f)\,\Pi(df)}.$$

Moreover, set $R_k(f) = \prod_{i=1}^k f(X_i)/f_0(X_i)$ for every $k = 1, 2, \ldots$. By exploiting the same martingale introduced in the proof of Theorem 4 in [14], one can show that the following holds.

PROPOSITION 1. *Suppose that* $K_{\varepsilon_n} < +\infty$ *and that*

(1) $$\sum_{n \geq 1} e^{-n\varepsilon_n^2/8} K_{\varepsilon_n} < +\infty,$$

*where* $(\varepsilon_n)_{n \geq 1}$ *is a sequence such that* $\varepsilon_n \to 0$ *and* $n\varepsilon_n^2 \to +\infty$. *Then*

(2) $$F_0^\infty\left(\liminf_n \left\{\sum_{j \geq 1} L_{n,j}^{(n)} < e^{-n\varepsilon_n^2/16}\right\}\right) = 1.$$

PROOF. One can easily check that the identity

$$L_{k+1,j}^{(n)}/L_{k,j}^{(n)} = \sqrt{f_{k,A_{n,j}}(X_{k+1})/f_0(X_{k+1})}$$

holds, given that

$$f_{k,A_{n,j}}(x) = \int_{A_{n,j}} f(x) R_k(f) \Pi(df) \Big/ \int_{A_{n,j}} R_k(f) \Pi(df), \qquad k \geq 1,$$

represents the predictive distribution restricted to the set $A_{n,j}$, whereas $f_{0,A_{n,j}}$ is the marginal density of the single observation restricted to $A_{n,j}$. Let $\mathscr{F}_k$ be the $\sigma$-algebra generated by the observations $X_1, \ldots, X_k$ and note that

$$E(L_{k+1,j}^{(n)}|\mathscr{F}_k) = L_{k,j}^{(n)}\{1 - h^2(f_{k,A_{n,j}}, f_0)/2\}.$$

Since $h(f_{k,A_{n,j}}, f_j) \leq \delta_n$, where $f_j$ is any density in $A_{n,j}$, from the triangle inequality one has $h(f_{k,A_{n,j}}, f_0) \geq \varepsilon_n - \delta_n = \gamma_n > 0$. Hence, we fix $k = n$ to obtain

$$E(L_{n+1,j}^{(n)}) \leq \sqrt{\Pi(A_{n,j})}(1 - \gamma_n^2/2)^{n+1}.$$

Choose a sequence $(\eta_n)_{n \geq 1}$ such that $\eta_n \to 0$ and $n\eta_n \to +\infty$. Apply Markov's inequality and the monotone convergence theorem to obtain

$$F_0^\infty\left(\sum_{j \geq 1} L_{n,j}^{(n)} > e^{-n\eta_n}\right) \leq e^{n\eta_n} E\left(\sum_{j \geq 1} L_{n,j}^{(n)}\right) \leq e^{n\eta_n} \sum_{j \geq 1} E(L_{n,j}^{(n)}).$$

Then

$$F_0^\infty\left(\sum_{j\geq 1} L_{n,j}^{(n)} > e^{-n\eta_n}\right) \leq \exp\{-n\{-\log(1-\gamma_n^2/2) - \eta_n\}\} \sum_{j\geq 1} \sqrt{\Pi(A_{n,j})}$$

$$= \exp\{-n\{-\log(1-\gamma_n^2/2) - \eta_n\}\} K_{\varepsilon_n}.$$

Setting $\phi_n = \varepsilon_n/2$ and $\eta_n = \varepsilon_n^2/16$, we have

$$-\log(1-\gamma_n^2/2) - \eta_n \geq \gamma_n^2/2 - \eta_n = \varepsilon_n^2/16.$$

Finally, condition (1) yields, by a straightforward application of the Borel–Cantelli lemma, the result in (2). □

Here we discuss a suitable lower bound for the denominator of the posterior, that is, $I_n := \int_{\mathbb{F}} R_n(f)\Pi(df)$. Regarding this point, previous contributions provide bounds in probability rather than almost surely. Indeed, Shen and Wasserman [13] and Ghosal, Ghosh and van der Vaart [6] give results of the type

$$I_n \geq \exp(-cn\varepsilon_n^2) \qquad \text{in } F_0^\infty\text{-probability}$$

for a constant $c > 0$, provided that $\Pi$ puts sufficient mass near $f_0$, where closeness is measured through a combination of the Kullback–Leibler divergence and the $L_2(F_0)$-norm of $\log(f_0/f)$. If $K(f, f_0) = \int \log(f_0(x)/f(x))f_0(x)\lambda(dx)$ and $V(f, f_0) = \int \{\log(f_0(x)/f(x))\}^2 f_0(x)\lambda(dx)$, then a neighborhood of the type above is defined as

$$(3) \qquad B(\varepsilon, f_0) = \{f : K(f, f_0) \leq \varepsilon^2, \ V(f, f_0) \leq \varepsilon^2\}.$$

One can now prove the following result.

THEOREM 1. *Suppose that $\varepsilon_n, \delta_n \to 0$ and $n\varepsilon_n^2, n\delta_n^2 \to +\infty$, and*

(i) $e^{-n\delta_n^2/16} K_{\delta_n} \to 0$;
(ii) *for some $C > 0$, $\Pi\{B(\varepsilon_n, f_0)\} \geq \exp(-Cn\varepsilon_n^2)$.*

*Then $\Pi_n(A_{\varepsilon_n}) \to 0$ in $F_0^\infty$-probability when $\delta_n \leq \phi\varepsilon_n$ for some sufficiently small $\phi > 0$.*

PROOF. Now,

$$\Pi_n(A_{\varepsilon_n}) \leq \sum_{j\geq 1} \Pi_n(A_{n,j}) \leq \sum_{j\geq 1} \sqrt{\Pi_n(A_{n,j})} = \sum_{j\geq 1} L_{n,j}^{(n)}/\sqrt{I_n}$$

and so

$$\Pi_n(A_{\varepsilon_n}) \leq \exp[-n\{\varepsilon_n^2/16 + n^{-1}(\log I_n)/2\}]$$

in $F_0^\infty$-probability. Moreover, by Lemma 8.1 in [6], condition (ii) implies that

$$I_n \geq \exp\{-n(1+C)\varepsilon_n^2\} \qquad \text{in } F_0^\infty\text{-probability.}$$

Hence, $n(\varepsilon_n^2/8 + n^{-1}\log I_n) \to +\infty$ in $F_0^\infty$-probability when $\delta_n \leq \phi\varepsilon_n$ and $(1+C)\phi^2 < 1/8$. The result follows. $\square$

A sequence $(\varepsilon_n)_{n\geq 1}$ satisfying (ii) in Theorem 1 is also referred to as *prior concentration* rate. As a simple illustration of condition (i), one can consider the discrete prior which puts mass $\Pi_k$ on the density $f_k$. If $\sum_{k\geq 1}\sqrt{\Pi_k} < +\infty$ then $K_{\varepsilon_n}$ is bounded by this sum and hence condition (i) reduces to $e^{-n\varepsilon_n^2/16} \to 0$, which is trivially satisfied for $\varepsilon_n = \lambda_n/\sqrt{n}$ for any $\lambda_n \to \infty$.

Before moving on to consider specific priors, we need to modify the above results, relying on the technique of Lijoi, Prünster and Walker [9], which was developed for establishing consistency of the mixture of Dirichlet process model. Let $N(\delta, \mathcal{S}, d)$ denote the minimum number of balls of radius at most $\delta$, with respect to the metric $d$, needed to cover the space $\mathcal{S}$. This is also known as the $\delta$-*covering number* of $\mathcal{S}$. Moreover, introduce a collection of sets $\{B_{n,k}: k\geq 1\}$ which, for any $n\geq 1$, forms a partition of $\mathbb{F}$. Accordingly, we denote by $\{A_{n,k,j}: j=1,\ldots,N(\varepsilon_n, B_{n,k}, h)\}$ an $\varepsilon_n$-covering of $B_{n,k}$ with respect to the Hellinger distance $h$. Hence, one can easily check that

$$\begin{aligned}\sqrt{\Pi(B_{n,k})} &= \sqrt{\sum_{j=1}^{N(\varepsilon_n, B_{n,k}, h)} \Pi(A_{n,k,j})} \\ &\geq \frac{1}{N(\varepsilon_n, B_{n,k}, h)} \sum_{j=1}^{N(\varepsilon_n, B_{n,k}, h)} \sqrt{\Pi(A_{n,k,j})}.\end{aligned} \qquad (4)$$

Next, it is clear that the family $\{A_{n,k,j}: j=1,\ldots,N(\varepsilon_n, B_{n,k}, h), k\geq 1\}$ is a partition of $\mathbb{F}$ into sets of diameter, with respect to the Hellinger distance, at most $\varepsilon_n$. Finally, using (4), one can write

$$K_{\varepsilon_n} = \sum_{k=1}^{\infty}\sum_{j=1}^{N(\varepsilon_n, B_{n,k}, h)} \sqrt{\Pi(A_{n,k,j})} \leq \sum_{k=1}^{\infty} N(\varepsilon_n, B_{n,k}, h)\sqrt{\Pi(B_{n,k})}.$$

Hence, we are interested in establishing, for some sequence $(\varepsilon_n)_{n\geq 1}$ such that $\varepsilon_n \to 0$ and $n\varepsilon_n^2 \to +\infty$ as $n\to\infty$, the validity of

$$e^{-n\varepsilon_n^2}\sum_{k=1}^{\infty} N(\varepsilon_n, B_{n,k}, h)\sqrt{\Pi(B_{n,k})} \to 0.$$

**3. Illustrations.** In the examples that follow, we show that the rate of convergence is governed by the concentration rate. In particular, we look at mixtures of Dirichlet processes and the random Bernstein polynomial model.

3.1. *Normal mixture of Dirichlet process.* The most widely used prior distribution for density estimation is undoubtedly the normal mixture of Dirichlet process (MDP) introduced by Lo [10] and later popularized by Escobar and West [3]. Such a random density function is given by

$$f_{\sigma,P}(x) = \phi_\sigma * P = \int \phi_\sigma(x-\theta) P(d\theta), \tag{5}$$

where the kernel $\phi_\sigma$ is the density function of the normal distribution with mean zero and variance $\sigma^2$. Moreover, $P$ is a Dirichlet process with parameter measure $\alpha(\cdot)$, and $\sigma$ has a prior distribution which we denote by $\mu$. The issue of strong consistency for the model (5) has been studied in [5] and [9], whereas rates are determined in [7].

We focus on the case in which the support of $\mu$ coincides with the interval $[\underline{\sigma}, \overline{\sigma}]$, where $0 < \underline{\sigma} < \overline{\sigma} < +\infty$, and we suppose that $f_0 = \phi_{\sigma_0} * P_0$. This is the same setting considered in [7]: when either $P_0$ has compact support or $\alpha$ has sub-Gaussian tails, they achieve the best rate of $(\log n)^\kappa/\sqrt{n}$ for $\kappa \geq 1$. In particular, for the usual Gaussian tails for $\alpha$, $\kappa = 3/2$. Although these models allow the desirable prior concentration rate [condition (ii) of Theorem 1] of $(\log n)/\sqrt{n}$, the worse entropy estimate determines their rate.

We, on the other hand, can obtain the target rate $(\log n)/\sqrt{n}$ in more general models. To this end, we introduce sets of the type

$$\mathscr{F}_{\underline{\sigma},a,\delta}^{\overline{\sigma}} = \bigcup_{\underline{\sigma} \leq \sigma \leq \overline{\sigma}} \{\phi_\sigma * P : P([-a,a]) \geq 1 - \delta\},$$

where $a > 0$. Moreover, recall that the inequality $h^2(f,g) \leq \|f-g\|_1$ yields

$$N(\sqrt{\delta}, \mathscr{F}^*, h) \leq N(\delta, \mathscr{F}^*, \|\cdot\|_1) \tag{6}$$

for any collection of density functions $\mathscr{F}^*$. Now, from [5] the upper bound for the $L_1$-metric entropy of set $\mathscr{F}_{\underline{\sigma},a,\delta}^{\overline{\sigma}}$ is given by

$$\log N(\delta, \mathscr{F}_{\underline{\sigma},a,\delta}^{\overline{\sigma}}, \|\cdot\|_1) \leq aC_\delta,$$

where $C_\delta = K\underline{\sigma}^{-1}\delta^{-1}\log(1/\delta)$ for some constant $K$. Hence, in view of (6), one finds that

$$\log N(\delta, \mathscr{F}_{\underline{\sigma},a,\delta^2}^{\overline{\sigma}}, h) \leq -\frac{2aK}{\underline{\sigma}} \delta^{-2} \log(\delta). \tag{7}$$

Now, for each $n$ let $(a_{n,j})_{j \geq 1}$ be an increasing sequence of positive numbers such that $a_{n,j} \uparrow +\infty$ as $j \to +\infty$, and for $j \geq 2$ set

$$\mathscr{G}_{\underline{\sigma},a_{n,j},\delta_n^2}^{\overline{\sigma}} := \bigcup_{\underline{\sigma} \leq \sigma \leq \overline{\sigma}} \{\phi_\sigma * P : P([-a_{n,j+1}, a_{n,j+1}]) \geq 1 - \delta_n^2,$$

$$P([-a_{n,j}, a_{n,j}]) < 1 - \delta_n^2\},$$

while setting
$$\mathscr{G}^{\overline{\sigma}}_{\underline{\sigma},a_{n,1},\delta_n^2} := \bigcup_{\underline{\sigma} \leq \sigma \leq \overline{\sigma}} \{\phi_\sigma * P : P([-a_{n,1},a_{n,1}]) > 1 - \delta_n^2\}.$$

Such sets cover the support of the distribution of the mixture of Dirichlet process defined in (5). It is obvious that for $j \geq 2$, $\mathscr{G}^{\overline{\sigma}}_{\underline{\sigma},a_{n,j},\delta_n^2}$ is included in $\mathscr{F}^{\overline{\sigma}}_{\underline{\sigma},a_{n,j},\delta_n^2}$, so

$$\log N(\delta_n, \mathscr{G}^{\overline{\sigma}}_{\underline{\sigma},a_{n,j},\delta_n^2}, h) \leq C_{\delta_n^2} a_{n,j}.$$

This suggests that for each $j \geq 2$, $\mathscr{G}^{\overline{\sigma}}_{\underline{\sigma},a_{n,j},\delta_n^2}$ has a finite Hellinger $\delta_n$-covering $\{C_{n,j,l} : l = 1, \ldots, N_{n,j}\}$, where $N_{n,j} \leq [\exp(C_{\delta_n^2} a_{n,j})] + 1$ and $[x]$ stands for the integer part of $x > 0$. Hence, setting

$$B_{n,j} = \{P : P([-a_{n,j+1},a_{n,j+1}]) \geq 1 - \delta_n^2, P([-a_{n,j},a_{n,j}]) < 1 - \delta_n^2\}$$

for $j \geq 2$, one has

$$K_{\delta_n} \leq N_{n,1} + \sum_{j \geq 2} N_{n,j} \sqrt{\Pi(\mathscr{G}^{\overline{\sigma}}_{\underline{\sigma},a_{n,j},\delta_n^2})} \leq N_{n,1} + \sum_{j \geq 2} N_{n,j} \sqrt{\mathscr{D}_\alpha(B_{n,j})},$$

where $\mathscr{D}_\alpha$ is the law of the Dirichlet process with parameter $\alpha$. If $V_{n,j} := [-a_{n,j},a_{n,j}]^c$, then $B_{n,j} \subset \{P : P(V_{n,j}) > \delta_n^2\}$. By the Markov inequality,

$$\mathscr{D}_\alpha(B_{n,j}) \leq \mathscr{D}_\alpha(\{P : P(V_{n,j}) > \delta_n^2\}) \leq c\alpha(V_{n,j})/\delta_n^2$$

for some constant $c$ and, thus, one has

$$K_{\delta_n} \leq N_{n,1} + \sum_{j \geq 2} c^{1/2} \delta_n^{-1} \exp\{a_{n,j} C_{\delta_n^2} - a_{n,j}^2\}$$

when, as we assume, $\alpha([-a,a]^c) \leq \exp(-2a^2)$. Hence, if we put, for $j \geq 2$, $a_{n,j} = jC_{\delta_n^2}$, then the summand in the bound for $K_{\delta_n}$ is bounded by

$$c^{1/2} \sum_{j \geq 2} \delta_n^{-1/2} \exp\{-(j-1)^2 C_{\delta_n^2}^2\},$$

which goes to zero as $n \to +\infty$. Now, $N_{n,1}$ is the Hellinger $\delta_n$-covering number of the set $\{P : P([-a_{n,1},a_{n,1}]) > 1 - \delta_n^2\}$. According to [7], it is the case that

$$N_{n,1} \leq L_1 \exp\{(\log(1/\delta_n))^2\}$$

when $a_{n,1} \leq L_2 \sqrt{(\log(1/\delta_n))}$ for constants $L_1$ and $L_2$. So, $\exp\{-n\delta_n^2\} K_{\delta_n} \to 0$ when $\exp\{-n\delta_n^2 + (\log(1/\delta_n))^2\} \to 0$, which occurs when $n\delta_n^2 - (\log(1/\delta_n))^2 \to +\infty$. Clearly, $\delta_n = M(\log n)/\sqrt{n}$ for some large enough $M$ is sufficient. Hence, for example, we obtain an overall rate of convergence $(\log n)/\sqrt{n}$ with normal $\alpha$ when the true mixing distribution $P_0$ has sub-Gaussian tails. This improves on Ghosal and van der Vaart [7], who obtain a rate of $(\log n)^{3/2}/\sqrt{n}$ in this case.

3.2. *Random Bernstein polynomials.* Another important prior for density estimation is the so-called *random Bernstein polynomial* introduced in [11]. Such a random density admits the representation

$$b(x;k,F) = \sum_{j=0}^{k}[F(j/k) - F((j-1)/k)]\beta(x;j,k-j+1),$$

where $\beta(x;a,b)$ is the beta density with parameters $a,b>0$. In the previous representation, $F$ is a random distribution function, usually chosen to be a Dirichlet process, and $k$ has distribution $p$ and is independent of $F$. Assuming $f_0$ is in the Kullback–Leibler support of the prior, consistency of such priors has been established in [12] and [14], where it has been shown that strong consistency holds under a suitable tail condition on $p$. Rates of convergence have been determined in [4]. In Theorem 2.3 of [4], it is shown that the prior concentration rate is $(\log n)^{1/3}/n^{1/3}$ and that the entropy rate is $(\log n)^{5/6}/n^{1/3}$, thus leading to an overall convergence rate of $(\log n)^{5/6}/n^{1/3}$.

Following our bound for $K_{\varepsilon_n}$ in Section 2, define $B_j$ to be the set of Bernstein polynomials of order $j$. Using the upper bound $N(\varepsilon_n, B_j, h) \leq (C/\varepsilon_n)^j$ provided by [4], we have

$$K_{\varepsilon_n} \leq \sum_{r=1}^{a_{n,1}}(C/\varepsilon_n)^r\sqrt{p_r} + \sum_{j=1}^{\infty}\sum_{r=a_{n,j}+1}^{a_{n,j+1}}(C/\varepsilon_n)^r\sqrt{p_r}.$$

Here we have introduced, for each $n$, an increasing sequence of reals $(a_{n,j})_{j\geq 1}$ which will be determined later on. Using the inequality $\sum_{r=1}^{M}c^r \leq Mc^M$ for $c>1$, we have

$$K_{\varepsilon_n} \leq a_{n,1}(C/\varepsilon_n)^{a_{n,1}} + \sum_{j=1}^{\infty}(C/\varepsilon_n)^{a_{n,j+1}}(a_{n,j+1}-a_{n,j})\sqrt{p_{a_{n,j}}}.$$

Here we have assumed that the $p_k$'s are decreasing for all large $k$ and we will also assume that $p_k \leq \exp(-4k\log k)$ for all large $k$. Therefore, putting $a_{n,j} = Cj/\varepsilon_n$, we have the summand term for the bound of $K_{\varepsilon_n}$ bounded by

$$C\varepsilon_n^{-1}\sum_{j=1}^{\infty}\exp\{(j+1)C\varepsilon_n^{-1}\log(C/\varepsilon_n) - 2jC\varepsilon_n^{-1}\log(jC/\varepsilon_n)\},$$

which is bounded by $C\varepsilon_n^{-1}\sum_{j=1}^{\infty}\exp\{-2jC\varepsilon_n^{-1}\log j\}$. In turn, this sum is bounded by $D/\varepsilon_n$ as $n \to +\infty$ for some constant $D$; the term $j=1$ ensures this. Returning to the first term in the bound for $K_{\varepsilon_n}$, we are interested in finding $\varepsilon_n$ for which $\varepsilon_n^{-1}\exp\{-n\varepsilon_n^2\} \to 0$ and $\exp\{-n\varepsilon_n^2 + C\varepsilon_n^{-1}\log(C/\varepsilon_n)\} \to 0$, when $n\varepsilon_n^2 - C\varepsilon_n^{-1}\log(C/\varepsilon_n) \to +\infty$. This clearly happens when

$$\varepsilon_n = M(\log n)^{1/3}/n^{1/3}$$

for sufficiently large $M$. Consequently, under the conditions of Theorem 2.3 in [4], we obtain a rate of convergence of $(\log n)^{1/3}/n^{1/3}$, which is the rate of convergence for the sieve MLE, whereas Ghosal [4] obtains a rate of $(\log n)^{5/6}/n^{1/3}$. Note that with lighter tails for $p$, namely $p_k < \exp(-2k^2)$, we can obtain a rate of $(\log n)/\sqrt{n}$ for $\varepsilon_n$, but the overall rate will remain at $(\log n)^{1/3}/n^{1/3}$.

**Acknowledgments.** We wish to thank Catia Scricciolo for helpful remarks. We are grateful to two referees for comments that led to improvements and, in particular, to one referee who pointed out an inaccuracy in an earlier version of the paper.

S. G. Walker
Institute of Mathematics, Statistics
  and Actuarial Science
University of Kent
Kent CT2 7NZ
United Kingdom
E-mail: S.G.Walker@kent.ac.uk

A. Lijoi
Dipartimento di Economia Politica
  e Metodi Quantitativi
Università degli Studi di Pavia
Via San Felice 5
27100 Pavia
Italy
E-mail: lijoi@unipv.it

I. Prünster
Dipartimento di Statistica e Matematica Applicata
Università degli Studi di Torino
Piazza Arbarello 8
10122 Torino
Italy
E-mail: igor@econ.unito.it